\documentclass{amsart}

\usepackage[cmtip,all]{xy}

\usepackage{amsmath,amssymb,amsthm, multicol,amsfonts,mathtools,stmaryrd, longtable,tabularray}
\usepackage{float} %Allows use of H (exact location) for figures
\usepackage{tikz,xcolor,soul,csquotes}
\usepackage[outline]{contour}
\contourlength{1.2pt}
\usetikzlibrary{cd,shapes,arrows,decorations.markings}
\usepackage{adjustbox} %to rescale tikzcd diagrams
\usepackage{paralist}
\usepackage[shortlabels]{enumitem}
\usepackage[textsize=footnotesize]{todonotes}
\usepackage[hidelinks]{hyperref}

\theoremstyle{theorem}
\newtheorem{thm}{Theorem}[section]
\newtheorem*{thm*}{Theorem}
\newtheorem{lem}[thm]{Lemma}
\newtheorem*{lem*}{Lemma}
\newtheorem{cor}[thm]{Corollary}
\newtheorem{prop}[thm]{Proposition}
\newtheorem*{prop*}{Proposition}
\newtheorem*{conj}{Conjecture}

\theoremstyle{definition}

\newtheorem{example}[thm]{Example}

\numberwithin{equation}{section}

\newcommand{\C}{\mathbb{C}}
\newcommand{\R}{\mathbb{R}}

\newcommand{\Z}{\mathbb{Z}}

\newcommand{\G}{\mathcal{G}}
\newcommand{\into}{\hookrightarrow}

\renewcommand{\subset}{\subseteq}
\newcommand{\set}[1]{\left\{#1\right\}}
\renewcommand{\st}{\,\middle\vert\,}
\newcommand{\parens}[1]{\left(#1\right)}
\newcommand{\bracks}[1]{\left[#1\right]}

\newcommand{\abs}[1]{\left|#1\right|}
\newcommand{\iso}{\cong}
\newcommand{\im}{\operatorname{im}}

\newcommand{\rank}{\operatorname{rank}}

\newcommand{\dad}{\operatorname{dad}}

\newcommand{\Vect}{\operatorname{Vect}}
\newcommand{\diag}{\operatorname{diag}}
\newcommand{\Tr}{\operatorname{Tr}}

\let\iso\cong
\let\cong\equiv
\let\phi=\varphi

\begin{document}

% \title[short text for running head]{full title}
\title{A Low-Dimensional Counterexample to the HK-Conjecture}

%    Only \author and \address are required; other information is
%    optional.  Remove any unused author tags.

%    author one information
% \author[short version for running head]{name for top of paper}
\author{Rachel Chaiser}
\address{Department of Mathematics
Kalamazoo College
1200 Academy Street,
Kalamazoo, MI 49006-3295,
USA}
%\curraddr{}
\email{rachel.chaiser@kzoo.edu}
\thanks{This work was partially supported by NSF Grants 2000057 and 2247424.}

%    \subjclass is required.
\subjclass[2020]{46L80, 22A22}

%\date{}

%\dedicatory{}

%    "Communicated by" -- provide editor's name; required.
%\commby{}

\begin{abstract}
We provide a counterexample to the HK-conjecture using the flat manifold odometers constructed by Deeley. Deeley's counterexample uses an odometer built from a flat manifold of dimension 9 and an expansive self-cover. We strengthen this result by showing that for each dimension $d\geq 4$ there is a counterexample to the HK-conjecture built from a flat manifold of dimension $d$. Moreover, we show that this dimension is minimal, as if $d\leq 3$ the HK-conjecture holds for the associated odometer. We also discuss implications for the stable and unstable groupoid of a Smale space.
\end{abstract}

\maketitle

\section{Introduction}

Matui's HK-conjecture \cite{Matui}, stated below, predicts a strong relationship between the homology of a certain class of groupoids and the $K$-theory of the reduced groupoid $C^*$-algebra. While there are a number of positive results \cite{dadhk, farsi, Matui, MR4170644, PY1, MR4052913}, the conjecture is known to fail in general.

\begin{conj}[HK-Conjecture]
Let $\G$ be a groupoid that is second countable, locally compact Hausdorff, \'etale, essentially principal, minimal, and ample. Then
\[
K_*(C^*_r(\G)) \iso \bigoplus_i H_{2i+*}(\G) \qquad \text{for $*=0,1$.}
\]
\end{conj}

In 2020, Scarparo \cite{Scarparo} provided the first counterexample to the HK-conjecture using a class of groupoids called odometers. While Scarparo's counterexample is essentially principal, in 2023 Deeley \cite{robin} gave the first counterexample in the principal case. Deeley's counterexample is an odometer built from dynamics on a flat manifold of dimension 9. We construct a counterexample to the HK-conjecture using the same class of flat manifold odometers but using a manifold of dimension 4. Moreover, we show that this dimension is minimal.

We begin with a flat manifold $Y$ and an expansive self-covering map $g:Y\to Y$. The existence of such a map is guaranteed by Epstein--Shub \cite{ES} and Shub \cite{Shub}. We then construct an action groupoid from the action of $\pi_1(Y)$ on the Cantor set
\[\Omega = \varprojlim\parens{\frac{\pi_1(Y)}{g_*^k\pi_1(Y)}, \text{ coset inclusion}}.
\]
The groupoid $\G=\Omega\rtimes\pi_1(Y)$ is called an odometer. It is principal and satisfies the hypotheses of the HK-conjecture.

Using results of \cite{reu,robin,Scarparo}, the homology of this odometer and the $K$-theory of the reduce groupoid $C^*$-algebra can be computed in terms of the corresponding invariants for $Y$ and the transfer map:
\[
H_*(\G) \iso \varinjlim(H_*(Y), \operatorname{tr}_H), \qquad \text{and}\qquad K_*(C^*_r(\G)) \iso \varinjlim(K_*(Y), \operatorname{tr}_K).
\]

One of the main results of the paper is Theorem \ref{thm:counterexample}, stated below. Throughout the paper, $T(G)$ denotes the torsion subgroup of a group $G$.

\begin{thm*}
Suppose that $Y$ is a flat manifold with $T(K^*(Y))\not\iso T\parens{\bigoplus_i H^{2i+*}(Y)}$. Then there exists an expansive self-cover $g:Y\to Y$ such that the corresponding odometer $\G$ is a counterexample to the HK-conjecture. Moreover, the relevant groupoid is principal.
\end{thm*}

To find such a flat manifold $Y$, we leverage the Atiyah--Hirzebruch spectral sequence which relates $K$-theory and cohomology. If $\dim(Y)\leq 4$, the spectral sequence collapses to two short exact sequences
\begin{align*}
& 0 \to H^4(Y) \to \widetilde{K}^0(Y) \to H^2(Y)\to 0,\\
& 0 \to H^3(Y) \to K^1(Y) \to H^1(Y)\to 0.
\end{align*}
The second will always split, while the first can be a nontrivial extension. Our key to finding a nontrivial extension is proved in Lemma \ref{lem:extension}. Our counterexample comes from the class of flat manifolds called real Bott manifolds, for which the cohomology ring $H^*(-;\Z_2)$ is readily computed using a binary matrix that defines the manifold. Notably, Deeley's counterexample also uses a real Bott manifold.

In Example \ref{xex1}, we present a real Bott manifold of dimension 4 that gives a counterexample to the HK-conjecture. In Theorem \ref{thm:largedim}, we extend this to counterexamples using a manifold of dimension $d$ for each $d\geq 4$. This is a strengthening of the results in \cite{robin}, where counterexamples are given in each dimension $d\geq 9$. In Theorem \ref{thm:dim3}, we prove that if $\dim(Y)\leq3$ then the HK-conjecture holds for the associated odometer.

In Lemma \ref{lem:dadodom}, we prove that the dynamic asymptotic dimension of our flat manifold odometers, denoted $\dad(\G)$, satisfies $\dad(\G)=\dim(Y)$. Thus, for each $d\geq 4$ there is counterexample to the HK-conjecture coming from a groupoid with $\dad(\G)=d$. This result is a companion to Theorem 4.19 of \cite{dadhk}, which says that the HK-conjecture holds for principal groupoids with $\dad(\G)\leq 2$ and $H_2(\G)$ free.

We conclude the paper by discussing implications of the main results to the groupoids associated to Smale spaces in Section \ref{sec:smale}.

This paper is based on the author's PhD thesis \cite{thesis}. We refer the reader to the thesis for examples and for more detail on the content of the paper.

\section*{Acknowledgments}
Many thanks to my PhD advisor, Robin Deeley, for his guidance and feedback throughout the project. Thanks to Rufus Willett and Agn\`es Beaudry for a thorough reading of my thesis and for their insightful comments. Thanks as well to the other members of my thesis committee. 

%The structure of the paper is as follows. In Section 2, we discuss flat manifolds, and the particular class of flat manifolds called real Bott manifolds from which we construct our counterexample. Section 3 contains an overview of the construction of odometers from flat manifolds, as well as their homology and $K$-theory. Section 3 also contains a discussion of dynamic asymptotic dimension. In Section 4 we discuss the necessary topological tools; these are the Atiyah--Hirzebruch spectral sequence and characteristic classes. Section 5 contains our main results regarding counterexamples. In Examples \ref{xex1} and \ref{xex2}, we give two manifolds in dimension 4 that produce a counterexample to the HK-conjecture. We then extend these examples to construct a counterexample in dimension $d$ for each $d\geq 4$, and we show that if $d\leq 3$ then the HK-conjecture holds. Finally, Section 6 is a departure from odometers to discuss the implications of our results to  the groupoids and $C^*$-algebras built from Smale spaces.

%-------------------------------------
\section{Flat Manifolds}
A {\bf flat manifold} $Y$ is a manifold of the form $\R^n/\pi$ where $\pi$ is a cocompact, torsionfree, discrete subgroup of $O(n)\ltimes \R^n$ that acts freely and properly on $\R^n$. This definition agrees with the definition that $Y$ be Riemannian flat (see \cite{charlap} Section II.5). The fundamental group of a flat manifold is a {\bf Bieberbach group} \cite{bieberbach}. Since the action of $\pi$ on $\R^n$ is free and proper, $\pi = \pi_1(Y)$ and that $Y$ is a model for $B\pi$. Every flat manifold $Y$ corresponds to a short exact sequence
\[0 \to \Z^n \to \pi \to F \to 0\]
where $F$ is a finite group called the {\bf holonomy}.% One obtains this short exact sequence by identifying the pure translations $\set{(I_n,p)\in\pi}\trianglelefteq \pi$ with the lattice $\Z^n$. Then $F$ is the set of rotations and reflections $\set{A\in O(n)\st (A,a)\in\pi}$.

It is a result of Epstein--Shub in \cite{ES} that every flat manifold admits a locally expansive self-covering map.
\begin{thm}
\label{thm:ES}
Let $Y$ be a flat manifold associated to the short exact sequence
\[0 \to \Z^n \to \pi \to F \to 0.\]
Then for any $k\in\Z$ such that $|k|F|+1|>1$ (i.e., $k>0$ or $k<-2/|F|$) the map $(k|F|+1)I_n:\R^n\to\R^n$ projects to an affine expanding endomorphism of $Y$.
\end{thm}
While this result was originally stated for $k>0$, their proof holds more generally for negative values of $k$ such that $|k|F|+1|>1$. So we choose to state the theorem in its full generality. By Lemma 1 of \cite{Shub}, such an expanding endomorphism is necessarily a covering map. Throughout, we call such a map an expansive self-cover and denote it by $g:Y \to Y$. See Examples 2.7 and 2.8 of \cite{reu} for explicit examples of expansive self-covers on flat manifolds.

\subsection{Real Bott Manifolds}

To construct our counterexample, we consider a class of flat manifolds called real Bott manifolds because the cohomology ring with $\Z_2$ coefficients is easily computed. Our overview follows \cite{masuda,KM}.

A {\bf real Bott tower} is a sequence of $\R P^1$ (i.e., $S^1$) bundles
\begin{center}
\begin{tikzcd}[column sep=2em]
M_n \arrow[r, "\R P^1"] & M_{n-1} \arrow[r, "\R P^1"] & \cdots \arrow[r, "\R P^1"] & M_2 \arrow[r, "\R P^1"] & S^1 \arrow[r, "\R P^1"] & pt
\end{tikzcd}
\end{center}
such that at each stage $M_j\to M_{j-1}$ is the projectivization of a real vector bundle of the form $L_{j-1}\oplus \underline{\R}$ over $M_{j-1}$, where $\underline{\R}$ is the trivial line bundle and $L_{j-1}$ is any real line bundle over $M_{j-1}$. The total space $M_n$ is a {\bf real Bott manifold}.

Since real line bundles are classified by their first Stiefel--Whitney class (this is discussed in Section \ref{sec:CC}), each $L_{j-1}$ corresponds to an element in $H^1(M_{j-1};\Z_2)\iso \Z_2^{j-1}$. Viewing these elements as column vectors, padding the bottoms with 0's to create a column vector in $\Z_2^n$, and lining them up in $j$-ascending order, we obtain an $n\times n$ upper triangular $\set{0,1}$-matrix that encodes the real Bott manifold. Such a matrix is called a {\bf Bott matrix}, and given a Bott matrix $A$ we can recover a real Bott manifold denoted $M(A)$ by reversing this construction. For examples, we refer the reader to \cite{thesis,KM}.

Much of the structure of $M(A)$ can be read off of the defining matrix $A$. The relevant results are stated below.

\begin{prop}[\cite{KM} Lemma 2.1]
\label{prop:rbmcohom}
As a graded ring, $H^*(M(A);\Z_2)$ is generated by degree one elements $x_1,\dots,x_n$ with $n$ relations
\[x_j^2 = x_j\sum_{i=1}^n A_{i,j}x_i \quad \text{for }j=1,\dots, n.\]
\end{prop}
%\todo{On p.4 of \cite{KM}, they say that $x_j =w_1(\gamma_j)$ where $\gamma_j$ is the ``canonical line bundle over $M_j$''.} %{\color{red} I don't know what the canonical line bundle over something that isn't projective space would be.}

\begin{prop}[\cite{KM} Lemma 2.2]
\label{prop:rbmorientable}
The real Bott manifold $M(A)$ is orientable if and only if for each row of $A$ the sum of the entries is $0\pmod{2}$.
\end{prop}

%We can use Proposition \ref{prop:rbmpi} to compute the first homology group of a real Bott manifold.
The following result may be well-known, but we include a proof for completeness.

\begin{lem}
\label{lem:H1}
Let $A$ be a Bott matrix. Let $r$ denote the number of columns of $A$ that have all zero entries, and let $t$ denote the number of columns of $A$ that have nonzero entries. Then
\[H_1(M(A);\Z) = \Z^r \oplus \Z_2^t.\]
Moreover, the  torsion subgroup of $H^2(M(A);\Z)$ is $\Z_2^t$.
\begin{proof}
Let $\pi = \pi_1(M(A))$. We will use that $H_1(M(A)) = \pi/[\pi,\pi]$. From Section 3 of \cite{KM}, the fundamental group $\pi\subset O(n)\ltimes \R^n$ is generated by the rigid motions $s_1,\dots,s_n$ defined by
\[s_i  = \parens{\operatorname{diag}\parens{ (-1)^{A_{i,1}},(-1)^{A_{i,2}},\dots, (-1)^{A_{i,n}}}, \frac{1}{2}e_i},\]
where $e_i$ is the $i$th standard basis vector for $\R^n$. %Denote $s_i = (R_i, \frac{1}{2}e_i)$. Recall that the group structure on $O(n)\rtimes \R^n$ is
%\[ (A,a)(B,b) = (AB,a+Ab),\qquad\text{and}\qquad (A,a)^{-1} = (A^{-1}, -A^{-1}a). \]
%Notice that $R_i^{-1}=R_i$ and $R_i^2 = I$ for each $i$, and that each $R_i$ and $R_j$ commute.
Without loss, assume $j>i$. A computation (see Lemma 5.2.7 in \cite{thesis}) shows
\[
(s_is_j)(s_i^{-1}s_j^{-1}) = \parens{I, \frac{1}{2}((-1)^{A_{ij}}e_j - e_j )},
\]
which is the identity if and only if $A_{ij}=0$. Moreover, if $A_{ij}\neq0$ so that $s_i$ and $s_j$ do not commute, then $s_is_js_i^{-1}s_j^{-1} = (I,-e_j)$ which implies that $s_is_js_i^{-1} = s_j^{-1}$. Upon abelianizing, this gives $s_j^2=1$ in $H_1(M(A))$. Thus, $s_j$ generates a factor of $\Z_2$ if and only if there is a 1 in column $j$, and it generates a factor of $\Z$ otherwise.

The moreover part follows from the Universal Coefficient Theorem.
\end{proof}
\end{lem}

%\begin{example}
%Consider the real Bott manifold $Y$ corresponding to the matrix
%\[\begin{bmatrix} 0 & 1 & 0 & 0 \\ 0 & 0 & 1 & 0 \\ 0 & 0 & 0 & 1 \\ 0 & 0 & 0 & 0 \end{bmatrix}.\]
%This is \#12 in the list of real Bott manifolds found in \cite{KM}. Applying Lemma \ref{lem:H1}, we see that $H_1(Y;\Z)=\Z\oplus\Z_2^3$; and by the UCT we have that $T(H^2(Y;\Z))=\Z_2^3$.
%\end{example}

\section{Odometers}
\label{sec:odom}

While odometers can be defined more generally (see for example \cite{Scarparo}), we are only concerned with odometers built from flat manifolds as in \cite{robin}. Start with $Y$ a flat manifold and $g:Y\to Y$ an expansive self-cover. There is a chain of subgroups $\pi_1(Y) \supset g_*\pi_1(Y) \supset g_*^2\pi_1(Y) \supset \cdots$ ordered by inclusion. Pulling back along the inclusion gives maps between quotient groups
\[
f^{k+1}_{k}: \frac{\pi_1(Y)}{g^{k+1}_*\pi_1(Y)} \to \frac{\pi_1(Y)}{g^{k}_*\pi_1(Y)}.
\]
The inverse limit
\[\Omega \coloneqq \varprojlim \parens{\frac{\pi_1(Y)}{g^k_*\pi_1(Y)},f^{k+1}_{k}}
\]
is a Cantor set since each quotient is finite and nontrivial. Then $\pi_1(Y)$ acts on $\Omega$ by group multiplication on the coset representatives. The action groupoid $\G = \Omega \rtimes\pi_1(Y)$ is called an {\bf odometer}.

Such an odometer $\G$ is principal and satisfies the hypotheses of the HK-conjecture. %That $\G$ is 2nd countable, locally compact Hausdorff, \'etale, and ample is clear. 
In particular, $\pi_1(Y)$ is always a discrete group so $\G$ is \'etale; and $\G^{(0)}=\Omega$ is a Cantor set by construction so that $\G$ is ample. That $\G$ is minimal is proven in \cite{CP}. That $\G$ is principal is proven in \cite{robin} Proposition 3.1, and relies on the fact that the covering map $g:Y\to Y$ is locally expansive.

\subsection{Homology and $K$-theory of Odometers}
\label{sec:odomH}
%Let $Y$ be a flat manifold, $g:Y\to Y$ an expansive self-cover, and $\G$ the associated odometer. 
For flat manifold odometers, both $H_*(\G)$ and $K_*(C^*_r(\G))$ can be computed from the inputted $Y$ and $g$. The following are synthesized from the results in \cite{reu,robin,Scarparo}.

\begin{prop}%[\cite{robin} Theorem 4.1]
\label{prop:odomH}
Let $\G$ be the odometer associated to a flat manifold $Y$ and expansive self-cover $g:Y\to Y$. Let $\operatorname{tr}_H$ (resp. $\operatorname{tr}_K$) denote the transfer map in homology (resp. $K$-homology). Then
\begin{align*}
H_*(\G)&\iso\varinjlim(H_*(Y),\operatorname{tr}_H)\\
K_*(C^*_r(\G))&\iso\varinjlim(K_*(Y),\operatorname{tr}_K).
\end{align*}
\end{prop}

Given a specific $n$-fold self-cover $g:Y\to Y$, $H_*(\G)$ is is straightforward to compute from $H_*(Y)$, since the transfer in homology satisfies $g_*\circ \operatorname{tr}_H=\times n$. In particular, see Example 4.4 of \cite{reu} for a concrete example of this type of computation. Crucial for the construction of our counterexample, if one uses $g$ as in Theorem \ref{thm:ES}, the torsion subgroup of $H_*(Y)$ is preserved in the inductive limit.

\begin{prop}[\cite{robin} Theorem 4.4]
\label{prop:torsionpreserved}
Let $Y$ be a flat manifold. Then there exists an expansive self-cover $g:Y \to Y$ such that
\[
T(\varinjlim(H_*(Y),\operatorname{tr}_H)) \iso T(H_*(Y)).
\]
\end{prop}

In contrast, the transfer in $K$-homology fails to satisfy $g_*\circ \operatorname{tr}_K = \times n$ in general (see \cite{transfer}), and so an explicit computation of $K_*(C^*_r(\G))$ is less accessible. However, the following lemma gives us enough information about $T(K_*(C^*_r(\G)))$ so that we do not need an explicit computation of $K_*(C^*_r(\G))$ for our counterexample.

\begin{lem}
\label{lem:torlim}
Let $G$ be a finitely generated abelian group and $\beta:G\to G$ a homomorphism. Then $T(\varinjlim(G,\beta))$ is isomorphic to a subgroup of $T(G)$. Moreover, if $\beta|_{T(G)}:T(G)\to T(G)$ is an isomorphism then $T(\varinjlim(G,\beta))$ is isomorphic to $T(G)$. 
\begin{proof}
Let $[\gamma,n]\in T(\varinjlim(G,\beta))$ be nonzero. Then there exists $k>1$ and $\ell>0$ with $k[\gamma,n]=[k\gamma,n] = [\beta^{\ell}(k\gamma),n+\ell]=[0,n+\ell]$. By assumption, $[\gamma,n] = [\beta^{\ell}(\gamma),n+\ell]\neq 0$, so $\beta^{\ell}(\gamma)\neq 0$ for any $\ell$. So we must have $k\gamma=0$ and hence $\gamma\in T(G)$.

By Proposition 2.5 of \cite{robin}, $T(\varinjlim(G,\beta))$ is a finite group. So we can list its elements as $[\gamma_1,n_1], [\gamma_2,n_2],\dots,[\gamma_m,n_m]$. By applying $\beta$, we can assume that $n_1 = n_2 = \cdots = n_m$. Call this common integer $n$. Then the map $\Phi: T(\varinjlim(G,\beta)) \to T(G)$ defined by $\Phi([\gamma_i,n])=\gamma_i$ is injective and hence $\im\Phi \iso T(\varinjlim(G,\beta))$.

For the moreover part, define $\Psi: T(G) \to T(\varinjlim(G,\beta))$ by $\gamma \mapsto [\gamma,0]$. Then $\Psi$ is injective since $\beta$ is, since $[\gamma,0]=0$ if only if $\beta^{\ell}(\gamma)=0$ for some $\ell>0$. For surjectivity, let $[\gamma,n]\in T(\varinjlim(G,\beta))$. Since $\beta$ is an isomorphism, we can form $\hat{\gamma} = \beta^{-k}(\gamma)$. Then $\Psi(\hat{\gamma}) = [\hat{\gamma},0] = [\gamma,k]$.
\end{proof}
\end{lem}

\subsection{Dynamic asymptotic dimension}
\label{sec:dad}
Dynamic asymptotic dimension is a dimension theory for groupoids and group actions. We refer the reader to \cite{dad} for the precise definition. Important for us is the fact, proven below, that for odometers built from flat manifolds, the dynamic asymptotic dimension coincides with the dimension of the inputted manifold. The relationship between dynamic asymptotic dimension and the HK-conjecture is covered in \cite{dadhk} and in Section \ref{sec:leq4} below. 

\begin{lem}
\label{lem:dadodom}
Let $\G$ be the odometer associated to a flat manifold $Y$ and an expansive self-cover $g:Y\to Y$ as in Theorem \ref{thm:ES}. Then $\operatorname{dad}(\G)=\dim(Y)$.
\begin{proof}
Set $d=\dim (Y)$. By Theorem 3.11 of \cite{reu}, $H_n(\G)=0$ for all $n>d$ and if $Y$ is orientable then $H_d(\G)=\Z$. If $Y$ is nonorientable, then $\Z_2\subset H_{d-1}(Y)$ and by Proposition \ref{prop:torsionpreserved} we have $\Z_2\subset H_{d-1}(\G)$ as well. Since either $H_d(\G)$ is free or $H_n(\G)$ is nontrivial for $n<d$, by Theorem 3.36 of \cite{dadhk} we have $\operatorname{dad}(\G)\geq d= \dim(Y)$.

For an upper bound on $\operatorname{dad}(\G)$, we use that $\G$ is Morita equivalent to the unstable groupoid $\G^u(P)$ of the Smale space $(X,\phi)$, where $X = \varprojlim(Y,g)$ and $\phi:X\to X$ is the shift map. This Morita equivalence is proven in Lemma 2.4 of \cite{reu} (see also Example 3.4 of \cite{PY3}), and stated explicitly in Section \ref{sec:smale} as Proposition \ref{prop:moritagroupoid}. In \cite{bonicke}, it is shown that dynamic asymptotic dimension is a Morita invariant. By Theorem 3.7 of \cite{DS}, $\operatorname{dad}(\G^u(P))\leq \dim(X)$. It is shown in \cite{coveringdim} that such an inverse limit satisfies $\dim(X) \leq \dim(Y)$, which gives the desired result. %follows from \url{https://www.jstor.org/stable/pdf/2044041.pdf} to get $\dim(X)\leq \dim(Y)$ and cohomological dimension gives $\dim(X)\geq \dim(Y)$.\todo{Do I need to include this proof for $\dim(X)=\dim(Y)$?}
\end{proof}
\end{lem}

\section{Topological Preliminaries}
In this section we collect a number of well-known topological results that are used for our main results. Throughout this section, $X$ denotes a  finite CW-complex.

\subsection{Atiyah--Hirzebruch Spectral Sequence}
\label{sec:user}
The Atiyah--Hirzebruch spectral sequence relates the $K$-theory of $X$ to its cohomology. The $E_2$-page of the spectral sequence is given by $E_2^{p,q}=H^p(X;K^q(pt))$, and if all differentials are trivial and there are no so-called extension issues, then one has $K^*(X) \iso \bigoplus_i H^{2i+*}(X)$. In \cite{robin}, a counterexample to the HK-conjecture is constructed by demonstrating nontrivial differentials. However, in dimensions at most 4 all differentials are trivial (see \cite{thesis} for an analysis of this) so we instead look for extension issues.

If $\dim(X)\leq 4$ the spectral sequence collapses to two short exact sequences
\begin{align*}
& 0 \to H^4(X) \to \widetilde{K}^0(X) \to H^2(X)\to 0,\\
& 0 \to H^3(X) \to K^1(X) \to H^1(X)\to 0.
\end{align*}
Since $H^1(X)$ is a free group, the second sequence splits and we have $K^1(X) \iso H^1(X)\oplus H^3(X)$. Notice that if $\dim(X)\leq 3$ then the first sequence reduces to an isomorphism and we have $K^0(X)\iso H^0(X)\oplus H^2(X)$. If $\dim(X) = 4$ however, the first sequence may be a nontrivial extension.

We also need the following relationship between the torsion subgroups of $K^*(X)$ and $\bigoplus_i H^{2i+*}(X)$. A proof of which can be found in Item (5) in Section 2.4 of \cite{robin}.
\begin{prop}
\label{prop:Ktorsion}
Let $X$ be a finite CW-complex. Then
\[
|T(K^*(X))| \leq \abs{T\parens{\bigoplus_i H^{2i+*}(X)}}.
\]
\end{prop}

\subsection{Characteristic Classes}
\label{sec:CC}

Characteristic classes are invariants of vector bundles defined by associating cohomology classes to them. Let $\Vect_{\R}^n(X)$ denote the set of rank $n$ real vector bundles over $X$, and $\Vect_{\C}^n(X)$ denote the set of rank $n$ complex vector bundles over $X$. If $n$ is omitted we mean all vector bundles regardless of rank. To a real vector bundle $\xi\in\Vect_{\R}^n(X)$ one associates the {\bf Stiefel--Whitney classes} $w_k(\xi)\in H^k(X;\Z_2)$. To a complex vector bundle $\eta\in\Vect_{\C}^n(X)$ one associates the {\bf Chern classes} $c_k(\eta)\in H^{2k}(X;\Z)$. These classes satisfy the following axioms, which we state for Stiefel--Whitney classes. (The axioms for Chern classes are formally the same.) %\todo{Probably shouldn't be a theorem, not sure how to format.}

\begin{thm}
There is a unique sequence of functions $(w_k)_{k\geq 0}:\Vect_{\R}(X) \to H^k(X;\Z_2)$ satisfying
\begin{compactenum}[(a)]
\item {\it (Naturality)} $w_k(f^*(\xi))=f^*(w_k(\xi))$,
\item {\it (Whitney sum formula)}
\[
w_k(\xi\oplus \zeta) = \sum_{i+j=k}w_i(\xi) w_j(\zeta), \qquad \text{where $w_0=1$,}
\]
\item {\it (Rank)} $w_k(\xi)=0$ for $k>\rank(\xi)$,
\item {\it (Nontriviality)} For the canonical line bundle $\xi\to \R P^1$, $w_1(\xi)$ is the generator of $H^1(\R P^1;\Z_2) = \Z_2$.
\end{compactenum}
\end{thm}

Below is the list of facts needed in the next section. The first three are easy consequences of the axioms. Let $n$ denote the trivial $n$-bundle (real or complex based on context).

\begin{enumerate}
    \item The characteristic classes of a trivial bundle are trivial.
    \item Characteristic classes are stable isomorphism invariants; i.e., $w_k(\xi\oplus n) = w_k(\xi)$. The corresponding statement holds for Chern classes.
    \item Let $\eta\in\Vect_{\C}^n(X)$. If $c_k(\eta)\neq 0$ for some $k>0$, then $\eta-n$ is nontrivial in $\widetilde{K}(X)$ since $\eta-n$ is not stably trivial.
    \item Real line bundles over $X$ are classified by their first Steifel--Whitney class in the sense that $w_1:(\Vect_{\R}^1(X),\otimes) \to H^1(X;\Z_2)$ is a group isomorphism. Similarly, complex line bundles over $X$ are classified by their first Chern class in the sense that $c_1:(\Vect_{\C}^1(X),\otimes) \to H^2(X;\Z)$ is a group isomorphism. See Proposition 3.10 of \cite{hatchervb} for a proof.
\end{enumerate}

These cohomology groups are related by the coefficient sequence
\[0\to\Z\xrightarrow{\times 2} \Z \xrightarrow{\rho}\Z_2\to 0,\]
which induces the long exact sequence in cohomology
\[
\cdots \to H^1(X;\Z) \xrightarrow{\times 2} H^1(X;\Z) \xrightarrow{\rho} H^1(X;\Z_2) \xrightarrow{\beta} H^2(X;\Z) \xrightarrow{\times 2} H^2(X;\Z) \to \cdots,
\]
where the connecting map $\beta$ is the Bockstein map.

\begin{enumerate}
    \setcounter{enumi}{4}
    \item For $\eta\in \Vect_{\C}^n(X)$, we can form a vector bundle $r\eta\in \Vect_{\R}^{2n}(X)$ called its {\bf realification} by identifying the fibers $\C^n$ with $\R^{2n}$. This real vector bundle satisfies $w_{2k+1}(r\eta)=0$, and $w_{2k}(r\eta)=\rho(c_k(\eta))$ (see Propostion 3.8 in \cite{hatchervb}).
    \item For $\xi\in \Vect_{\R}^n(X)$, we can form a vector bundle $c\xi\in\Vect_{\C}^n(X)$ called its {\bf complexification} by considering the fibers $\R^n$ inside of $\C^n$. Then the Bockstein map represents complexification of real line bundles in the sense that $c_1 (c\xi) = \beta \circ w_1 (\xi)$ (see \cite{bockref} p.19 item (8)).
\end{enumerate}

\section{The Counterexample}
\label{ch:xex}

\subsection{A Criterion for a Nontrivial Extension}

%It is worth noting that the proof of Lemma \ref{lem:extension} was inspired by, Adams' proof that $K^0(\R P^{2n}) \iso \Z\oplus \Z_{2^n}$ (see Theorem 7.3 of \cite{adams}). While the techniques in \cite{adams} were not adapted to our context, Adams' computation is nonetheless illuminating as a high-level strategy for detecting nontrivial extensions in the Atiyah--Hirzebruch spectral sequence.
Although our criterion for a nontrivial extension may be mysterious, we hope that the proof is illuminating, as well as the discussion at the beginning of Section \ref{sec:xex}.

\begin{lem}
\label{lem:extension}
Let $X$ be a nonorientable connected manifold with $\dim(X)=4$ such that $T(H^2(X;\Z)) = \Z_2^t$ for some $t>0$. If there exists $x\in H^1(X;\Z_2)$ such that $x^4\neq 0$ then $\Z_4\subset K^0(X)$ and $T(K^0(X))\not\iso T\parens{\bigoplus_i H^{2i}(X)}$.
\begin{proof}
Write $H^2(X;\Z)=\Z^{b_2}\oplus\Z_2^t$. The Atiyah-Hirzebruch spectral sequence for $X$ gives a short exact sequence of the form
\begin{equation}
\label{eq:ses}
0\to \Z_2 \to \widetilde{K}^0(X) \to \Z^{b_2}\oplus\Z_2^t \to 0.
\end{equation}
Thus, the torsion elements in $K^0(X)$ must have order dividing 4.% To show that $\Z_4\subset K^0(X)$, we will construct a torsion element $\nu$ in $K^0(X)$ with $2\nu\neq 0$.

By Item (6) in Section \ref{sec:CC}, the Bockstein map represents complexification of real line bundles. The long exact sequence defining $\beta$ looks like
\[
\cdots \to H^1(X;\Z) \xrightarrow{\times 2} H^1(X;\Z) \xrightarrow{\rho} H^1(X;\Z_2) \xrightarrow{\beta} \Z^{b_2}\oplus \Z_2^t \xrightarrow{\times 2} \Z^{b_2}\oplus\Z_2^t \to \cdots.
\]
By exactness, $\im(\beta) = \ker(\times 2)=\Z_2^t$, which we assume is nontrivial. Thus, for each $\xi\in\Vect_{\R}^1(X)$, its complexification satisfies $c_1(c\xi) = \beta\circ w_1(\xi)\in\Z_2^t$. Letting 1 denote the trivial complex line bundle, we form $c\xi-1\in\widetilde{K}^0(X)$. The short exact sequence (\ref{eq:ses}) above gives that $c\xi-1$ must be a torsion element.

Now, let $x\in H^1(X;\Z_2)$ satisfy $x^4\neq 0$. Let $\xi$ be the real line bundle over $X$ with $w_1(\xi)=x$. Set $\eta = c\xi$ and $\nu= \eta -1 = c\xi-1$. We show that $\nu$ has order 4. Let $r\eta$ denote the realification of $\eta$ to a rank 2 real bundle. By Item (5) in Section \ref{sec:CC}, $\rho(c_k(\eta)) = w_{2k}(r\eta) \in H^{2k}(X;\Z_2)$ and $w_{2k+1}(r\eta)=0$ for all $k$. From the Whitney sum formula,
\[c_2(2\eta) = c_2(\eta)+c_1(\eta)c_1(\eta) + c_2(\eta) = c_1(\eta)^2,\]
since $\eta$ has complex rank 1. Now, applying $\rho$, we see that
\[\rho(c_1(\eta)^2) = \rho(c_2(2\eta))= w_4(2r\eta) = \sum_{i+j=4}w_i(r\eta) w_j(r\eta) = w_2(r\eta)^2 = \rho(c_1(\eta))^2 ,\]
since $w_1(r\eta)=w_3(r\eta)=0$, and $w_4(r\eta)=0$ since $r\eta$ has real rank 2. Now, by Item (6) of Section \ref{sec:CC}, $c_1(\eta)=c_1(c\xi)=\beta\circ w_1(\xi)= \beta(x)$. It is well-known that $\rho\circ\beta = \operatorname{Sq}^1$, the first Steenrod square (see \cite{HatcherAT}). Hence $\rho \circ \beta(x)=x^2$. Thus, if $x^4=\rho(c_1(\eta))^2\neq 0$ then $c_2(2\eta)\neq 0$. So $2\eta$ is not stably trivial and hence $\nu$ has order 4 in $K$-theory.
\end{proof}
\end{lem}

\subsection{The Desired Manifold}
\label{sec:xex}

Let $Y$ be a real Bott manifold. By Lemma \ref{lem:H1}, $H^2(Y;\Z)=\Z^{b_2}\oplus \Z_2^t$ where $t\geq 0$. From the proof of Lemma \ref{lem:extension}, the complex line bundles over $Y$ that are obtained as the complexification of a real line bundle are classified by the $\Z_2^{t}$ term. Thus, if $t>0$, then there exist real line bundles over $Y$ with nontrivial complexification. Lemma \ref{lem:extension} allows us to use the ring structure of $H^*(Y;\Z_2)$ to detect which of these complexified real line bundles give order 4 elements in $K$-theory. We can compute the ring structure of $H^*(Y;\Z_2)$ from the Bott matrix defining $Y$ by using Proposition \ref{prop:rbmcohom}.

Before we give the counterexample, we need a theorem to ensure that a flat manifold satisfying Lemma \ref{lem:extension} actually lifts to a counterexample to the HK-conjecture. The following is a stronger, cohomological version of Corollary 4.6 of \cite{robin}.

\begin{thm}
\label{thm:counterexample}
Suppose that $Y$ is a flat manifold with $T(K^*(Y))\not\iso T\parens{\bigoplus_i H^{2i+*}(Y)}$. Then there exists an expansive self-cover $g:Y\to Y$ such that the corresponding odometer $\G$ is a counterexample to the HK-conjecture. Moreover, the relevant groupoid is principal.
\begin{proof}
Let $H_k(Y) = \Z^{b_k}\oplus T_{k}$ be a decomposition of $H_k(Y)$ into its free and torsion parts. The UCT gives $H^k(Y) \iso \Z^{b_k}\oplus T_{k-1}$. The UCT in $K$-theory gives an analogous statement but with indices taken modulo 2. In particular, if $T(K^0(Y))\not\iso T\parens{\bigoplus_i H^{2i}(Y)}$ then $T(K_{1}(Y))\not\iso T\parens{\bigoplus_i H_{2i+1}(Y)}$, and similarly for $K^1(Y)$.% Thus it suffices to consider cohomology and $K$-theory.

By Proposition \ref{prop:torsionpreserved}, an expansive self-cover $g:Y\to Y$ as in Theorem \ref{thm:ES} satisfies 
\[T(H_*(\G))\iso T(\varinjlim(H_*(Y),\operatorname{tr}_H)) \iso T(H_*(Y)).\]
Take $g$ to be such a self-cover. By Lemma \ref{lem:torlim}, $T(K_*(C^*_r(\G))) \iso T(\varinjlim(K_*(Y), \operatorname{tr}_K))$ is isomorphic to a subgroup of $T(K_*(Y))$. By Proposition \ref{prop:Ktorsion}, we have that $\abs{T(K^*(Y))} \leq \abs{T\parens{\bigoplus_i H^{2i+*}(Y)}}$. Thus, even if some torsion of $K^*(Y)$ is killed in the inductive limit, the assumption that $T(K^*(Y))\not\iso T\parens{\bigoplus_i H^{2i+*}(Y)}$ ensures that we also have $T(K_*(C^*_r(\G))) \not \iso T(H_*(\G))$. That the groupoid is principal is shown in \cite{robin} Proposition 3.1.
\end{proof}
\end{thm}

The following real Bott manifold satisfies Lemma \ref{lem:extension}, and hence gives a low-dimensional counterexample to the HK-conjecture.

\begin{example}
\label{xex1}
Let $Y$ denote the real Bott manifold associated to the matrix
\[\begin{bmatrix} 0 & 1 & 0 & 0 \\ 0 & 0 & 1 & 0 \\ 0 & 0 & 0 & 1 \\ 0 & 0 & 0 & 0 \end{bmatrix}.\]
This is \#12 in the list of 4-dimensional real Bott manifolds in \cite{KM}. Note that $Y$ is nonorientable by Proposition \ref{prop:rbmorientable}. From Lemma \ref{lem:H1}, we have $H^2(Y;\Z)=\Z^{b_2}\oplus \Z_2^3$ where $b_2$ denotes the second betti number. Thus, there is a short exact sequence
\[
0 \to \Z_2 \to \widetilde{K}^0(Y) \to \Z^{b_2+1}\oplus\Z_2^3 \to 0.
\]

By Proposition \ref{prop:rbmcohom}, the generators $x_1,x_2,x_3,x_4$ of $H^1(Y;\Z_2)$ satisfy the relations
\[x_1^2=0,\quad x_2^2 = x_2x_1,\quad x_3^2=x_3x_2,\quad x_4^2=x_4x_3.\]
Notice that
\[x_4^4=x_4^2x_3^2= x_4 x_3 x_3 x_2 = x_4 x_3 x_2 x_2 = x_4 x_3 x_2 x_1,\]
which is the generator of $H^4(Y;\Z_2)$. By Lemma \ref{lem:extension}, $\widetilde{K}^0(Y)$ is a nontrivial extension to the above sequence. Thus, we have
\begin{align*}
K^0(Y) &= \Z^{b_2+1}\oplus\Z_2^2\oplus\Z_4 \qquad \text{while}\qquad \bigoplus_{i}H^{2i}(Y) = \Z^{b_2+1}\oplus\Z_2^4.
%K^1(Y) &= \Z\oplus H^3(Y;\Z).
\end{align*}
Therefore, by Theorem \ref{thm:counterexample}, there exists an expansive self-cover $g:Y\to Y$ such that the associated odometer is a counterexample to the HK-conjecture.
\end{example}

While we omit the details, a similar analysis shows that the real Bott manifold associated to the matrix
\[\begin{bmatrix} 0 & 0 & 0 & 1 \\ 0 & 0 & 1 & 0 \\ 0 & 0 & 0 & 1 \\ 0 & 0 & 0 & 0 \end{bmatrix},\]
which is \#6 in the list of 4-dimensional real Bott manifolds in \cite{KM}, also gives a counterexample to HK-conjecture.

\subsection{Dimensions Greater than 4}
\label{sec:geq4}
Let $Y$ be the manifold from Example \ref{xex1}, and let $\mathbb{T}^n$ denote the $n$-torus. Then $Y\times \mathbb{T}^n$ gives a counterexample to the HK-conjecture in dimension $4+n$. Note that this is a strengthening of the main result of \cite{robin}, which gives a counterexample to the HK-conjecture in each dimension $d\geq 9$. %For brevity, let $H^{even}(Y)$ denote the direct sum of the even cohomology groups of $Y$, and similarly for $H^{odd}(Y)$.

\begin{thm}
\label{thm:largedim}
%For each $n\geq 1$ there is a flat manifold of dimension $4+n$ giving a counterexample to the HK-conjecture.
For each $d\geq 4$, there exists a principal groupoid $\G$ with $\dad(\G) =d$ such that $\G$ is a counterexample to the HK-conjecture.
\begin{proof}
Let $Y$ be the manifold from Example \ref{xex1} and let $n\geq 1$ be an integer. Using the K\"unneth Formula and the UCT, we compute
\begin{align*}
\bigoplus_i H^{2i}(Y\times \mathbb{T}^n)) &\iso \bigoplus_i H^i(Y)^{2^{n-1}},\\ % H^{even}(Y)^{2^{n-1}})\oplus T(H^{odd}(Y)^{2^{n-1}})\\
K^{0}(Y\times \mathbb{T}^n) &\iso K^{0}(Y)^{2^{n-1}}\oplus K^{1}(Y)^{2^{n-1}}.
\end{align*}
Hence $T(\bigoplus_i H^{2i}(Y\times \mathbb{T}^n)) \not\iso T(K^{0}(Y\times \mathbb{T}^n))$ since $T(\bigoplus_i H^{2i}(Y)) \not\iso T(K^{0}(Y))$. By Theorem \ref{thm:counterexample}, there is an odometer $\G$ built from $Y\times \mathbb{T}^n$ giving a counterexample to the HK-conjecture. By Lemma \ref{lem:dadodom}, $\dad(\G)=\dim(Y\times\mathbb{T}^n)=4+n$.
\end{proof}
\end{thm}

%\begin{cor}
%For each $d\geq 4$, there exists a principal groupoid $\G$ with $\dad(\G) =d$ such that $\G$ is a counterexample to the HK-conjecture.
%\end{cor}

\subsection{Dimensions Less than 4}
\label{sec:leq4}
The main result of this section is that if $\dim(Y)\leq 3$, then the HK-conjecture holds for the associated odometer $\G$. Thus 4 is the minimal dimension for a counterexample to the HK-conjecture coming from the flat manifold odometers considered in this paper.

First, we need some results on the structure of $H_*(\G)$. The following can be found as Theorems 3.10 and 3.11 in \cite{reu}.

\begin{lem}%[\cite{reu} Theorem 3.11]
\label{lem:HdFree}
Let $\G$ be the odometer associated to a flat manifold $Y$ of dimension $d$ and an expansive $n$-fold self-cover $g:Y\to Y$. Then
\begin{align*}
H_0(\G) &\iso \varinjlim(\Z,\times n) \iso \Z\bracks{\frac{1}{n}}\\
H_d(\G)&\iso
\begin{cases}
\Z & \text{if $Y$ is orientable,}\\
0 & \text{if $Y$ is nonorientable.}
\end{cases}
\end{align*}
\end{lem}

If $\dim(Y)\leq 3$, then $H_n(\G)=0$ for $n>3$ and the Proietti--Yamashita spectral sequence \cite{PY1} reduces to two short exact sequences
\begin{equation}
\label{eq:K0seq}
0\to H_0(\G) \xrightarrow{\mu_0} K_0(C^*_r(\G)) \to H_2(\G) \to 0  
\end{equation}
\begin{equation}
\label{eq:K1seq}
0\to H_1(\G) \xrightarrow{\mu_1} K_1(C^*_r(\G)) \to H_3(\G) \to 0.   
\end{equation}
In particular, see the proof of Theorem 4.19 in \cite{dadhk}. The maps $\mu_0$ and $\mu_1$ are defined in Lemmas 4.10 and 4.11 of \cite{dadhk}, respectively.

\begin{thm}
\label{thm:dim3}
Let $\G$ be the odometer associated to a flat manifold $Y$ and an expansive $n$-fold self-cover $g:Y\to Y$. If $\dim(Y)\leq 3$, then the HK-conjecture holds for $\G$. 
\begin{proof}
First, assume that $\dim(Y)\leq 2$. The result follows from Theorem 4.19 of \cite{dadhk} since $\dad(\G)=2$ by Lemma \ref{lem:dadodom}, and $H_2(\G)$ is free by Lemma \ref{lem:HdFree}. 

Now, assume that $\dim(Y) = 3$. Again, by Lemma \ref{lem:HdFree}, $H_3(\G)$ is either $\Z$ or 0. In either case, the sequence (\ref{eq:K1seq}) splits. To show that the sequence (\ref{eq:K0seq}) splits, we will construct a section of the map $\mu_0$.

By Proposition 4.14 of \cite{dadhk}, the map $\mu_0: H_0(\G) \to K_0(C^*_r(\G))$ is induced by the inclusion $\iota:C_0(\G^{(0)}) \into C^*_r(\G)$ as follows. Notice that $K_0(C_0(\G^{(0)})) \iso C(\G^{(0)},\Z)$. Then there is a natural projection $p: C(\G^{(0)},\Z) \to H_0(\G)$, 
%since $H_0(\G)\iso C(\G^{(0)},\Z)/\im d_1$, where $d_1$ is the first differential in homology. Thus,
and $\mu_0$ is the map making the following diagram commute.
\begin{center}
    \begin{tikzcd}
{C(\G^{(0)},\Z)} \arrow[d, "p"'] \arrow[rd, "\iota_*"] &                \\
H_0(\G) \arrow[r, "\mu_0"]                             & K_0(C^*_r(\G))
\end{tikzcd}
\end{center}
To construct a section of $\mu_0$, we will leverage the inductive limit structure on each of the groups in the commutative triangle, construct a section on the level of the constituent groups, then show that it commutes with the connecting maps and hence lifts to a section on the level of the inductive limit.

First, we establish some notation. Let $\pi=\pi_1(Y)$ and let $\Omega_k = \pi/g_*^k\pi$. Then $\G^{(0)}\iso\Omega = \varprojlim (\Omega_k, f_k^{k+1})$. Hence
\[
C(\G^{(0)},\Z) \iso \varinjlim (C(\Omega_k,\Z), -\circ f_{k}^{k+1}).
\]
Let $M_{n^k}=M_{n^k}(\C)$. Using Proposition 2.3 of \cite{Scarparo} (also see the proof of Theorem 4.2 in \cite{robin}),
\[
K_0(C^*_r(\G)) \iso \varinjlim (K_0(M_{n^k}\otimes C^*_r(\pi))), \phi_k ),
\]
where the connecting map is induced by the usual diagonal embedding $\phi_k:M_{n^k}\to M_{n^{k+1}}$. Lastly, from Lemma \ref{lem:HdFree} we have
\[
H_0(\G) \iso \varinjlim (\Z,\times n).
\]
We analyze the splitting for indices $k=1,2$ and show that it is given by the trace on $M_{n^k}\otimes C^*_r(\pi)$. The proof for general $k$ follows in a similar way.

Expanding the inductive limits, the commutative triangle defining $\mu_0$ induces a sequence of commutative triangles.
\begin{center}
\begin{tikzcd}
{C(\Omega_1,\Z)} \arrow[d, "p"] \arrow[dd, pos=0.4, "\iota_*", bend left=45] \arrow[r, "-\circ f_1^2"] & {C(\Omega_2,\Z)} \arrow[d, "p"] \arrow[r] \arrow[dd, pos=0.4, "\iota_*", bend left=45]   & \cdots \arrow[r] & {C(\Omega,\Z)} \arrow[d, "p"] \arrow[dd, "\iota_*", bend left=50] \\
\Z \arrow[d, "\mu_0"] \arrow[r, "\times n"]                                                  & \Z \arrow[d, "\mu_0"] \arrow[r]                                              & \cdots \arrow[r] & H_0(\G) \arrow[d, "\mu_0"]                                        \\
K_0(M_{n}\otimes C^*_r(\pi)) \arrow[r, "\phi_1"] \arrow[u, "\operatorname{Tr}", bend left]            & K_0(M_{n^2}\otimes C^*_r(\pi)) \arrow[r] \arrow[u, "\operatorname{Tr}", bend left] & \cdots \arrow[r] & K_0(C^*_r(\G)) \arrow[u, "\operatorname{Tr}", bend left]          
\end{tikzcd}
\end{center}
To understand the map $\mu_0$ and show that $\Tr$ is a section, we need to understand all of the other maps in the diagram. We will start at the top, with $C(\Omega_k,\Z)$. Since $\Omega_k$ is discrete, a basis for $C(\Omega_k,\Z)$ is given by indicator functions on the cosets. Now, $|\Omega_k|=n^k$ so we can enumerate the cosets via $\Omega_1 = \set{U_{i}}_{i=1}^n$ and $\Omega_{2} = \set{U_{ij}}_{i,j=1}^n$, where our convention will be that $f_1^2(U_{ij}) = U_i$. Then the connecting map in the sequence for $C(\Omega,\Z)$ is defined on basis elements by
\[
-\circ f_1^2: 1_{U_i} \mapsto \sum_{j=1}^n 1_{U_{ij}}.
\]

The vertical map $p: C(\Omega_k,\Z)\to \Z$ in the diagram is the component of the projection $C(\Omega,\Z) \to C(\Omega,\Z)/\im d_1 \iso H_0(\G)$, and it maps each basis element in $C(\Omega_k,\Z)$ to $1\in\Z$. It is clear that this makes the top left square of the diagram commute.

The inclusion map $\iota_*: C(\Omega_k,\Z) \to K_0(M_{n^k}\otimes C^*_r(\pi))$ is an embedding of $C(\Omega_k,\Z)$ as matrix units. Let $e_{11}\in M_{n^k}$ be a matrix unit and $u_e\in C^*_r(\pi)$ be the unit. Define $\iota_*$ by mapping each basis element in $C(\Omega_k,\Z)$ to the class of the projection $e_{11}\otimes u_e$, which has trace 1.

Now, the horizontal map $\phi_1$ is induced by $\phi_1(a) = \diag(a,\dots,a)$, where $a$ occurs $n$ times. Now, $\iota_*$ commutes with the connecting maps since any two projections $p\otimes u_e, q\otimes u_e \in M_{n^k}\otimes C^*_r(\pi)$ with $\Tr(p\otimes u_e)=\Tr(q\otimes u_e)$ are in the same $K$-theory class. To see this, notice that $\Tr(p\otimes u_e)=\Tr(q\otimes u_e)$ implies $\Tr(p)=\Tr(q)$. So $p$ and $q$ are Murray-von Neumann equivalent in $M_{n^k}$ and thus there is a matrix $v$ implementing the equivalence. The element $v\otimes u_e$ then implements a Murray-von Neumann equivalence between $p\otimes u_e$ and $q\otimes u_e $.% Thus, we have that 
%\begin{align*}
%\iota_*\circ(-\circ f_1^2) (1_{U_i}) &= \iota_*\parens{ \sum_{j=1}^n 1_{U_{ij}} } = \parens{\bigoplus_{k=1}^n e_{11}}\otimes u_e, \quad \text{and}\\
%\phi_1\circ\iota_* (1_{U_i}) &= \phi_1(e_{11}\otimes u_e) = \parens{\bigoplus_{k=1}^ne_{11}}\otimes u_e.   
%\end{align*}
%While both outputs are in different size matrix algebras, they have the same trace and hence are equivalent in $K$-theory.

Finally, we turn to the map $\mu_0: \Z \to K_0(M_{n^k}\otimes C^*_r(\pi))$. Since it is defined so that $\mu_0\circ p= \iota_*$, we can define $\mu_0(1) = e_{11}\otimes u_e$. The argument of the previous paragraph implies that this makes the lower left square of the diagram commute.

Thus, we can define a section of $\mu_0$ using the standard un-normalized trace
\[\Tr: K_0(M_{n^k}\otimes C^*_r(\pi)) \to \Z \quad\text{defined by}\quad \sum_{i,j=1}^{n^k} e_{ij}\otimes a_{ij} \mapsto \sum_{i=1}^{n^k}\Tr(a_{ii}),
\]
which is the sum of the coefficient of $u_e$ in each $a_{ii}\in C^*_r(\pi)$. This is indeed a section of $\mu_0$ since $\Tr \circ \mu_0(1) = \Tr (e_{11}\otimes u_e)=1$. Further, it commutes with the connecting maps since $\Tr(\phi_1(a)) = \Tr(\diag(a,\dots,a) = n\Tr(a)$. Thus, the trace lifts to a section of $\mu_0$ in sequence \ref{eq:K0seq}.
\end{proof}
\end{thm}

%The case $\dim(Y)=3$ is more subtle and must be handled separately. Recall that if $\dim(Y)\leq 3$ then $K^*(Y)\iso \bigoplus_i H^{2i+*}(Y)$ (see Section \ref{sec:user}). Since the transfer map in $K$-homology does not necessarily satisfy the equation $g_*\circ \operatorname{tr}_K=\times n$, where $g:Y\to Y$ is an $n$-fold cover, it is possible in general that this isomorphism could not be preserved when taking the inductive limit. However, will show that the isomorphism is in fact preserved.

%The analysis of the Proietti--Yamashita spectral sequence \cite{PY1} done in \cite{dadhk} implies that if $\G$ is a 2nd countable, ample, and principal with $\operatorname{dad}(\G)\leq 3$ there is are short exact sequences
%\[0\to H_0(\G) \xrightarrow{\mu_0} K_0(C^*_r(\G)) \to H_2(\G) \to 0\]
%(see the proof of Theorem 4.19 in \cite{dadhk}).  We show that if this short exact sequence splits for our flat manifold odometers, then the HK-conjecture also holds for them. Analysis of this short exact sequence will be completed in future work.

\section{The Stable and Unstable Algebra of a Wieler Solenoid}
\label{sec:smale}

In this section, we discuss the implications of the main results of this paper to Smale spaces. We will cover enough on Smale spaces so that our case is self-contained. For a general treatment, we refer the reader to \cite{reu, smalenotes}.

\subsection{Groupoids from Smale Spaces}

A Smale space consists of a compact metric space $X$ and a homeomorphism $\phi:X\to X$ such that each point $x\in X$ has a local product structure consisting of a set of points that get closer to $x$ as $\phi$ is iteratively applied, and a set of points that get closer to $x$ as $\phi^{-1}$ is iteratively applied. There are then equivalence relations on $X$:
\begin{enumerate}
\item {\bf stable equivalence:} $x\sim_s y$ if and only if $\displaystyle\lim_{n\to\infty} d(\phi^n(x),\phi^n(y))=0$
\item {\bf unstable equivalence:} $x\sim_u y$ if and only if $\displaystyle \lim_{n\to\infty} d(\phi^{-n}(x),\phi^{-n}(y))=0$.
\end{enumerate}
The stable equivalence class of a point $x\in X$ is denoted $X^s(x)$ and is called a {\bf stable set}; the unstable equivalence class is denoted $X^u(x)$ and called an {\bf unstable set}.

Let $\mathbf{P}$ be a finite set of $\phi$-invariant periodic points. If $\phi$ has a fixed point, which we have below, then we will take $\mathbf{P}=\set{P}$. Define
\begin{align*}
X^u(\mathbf{P}) &= \set{x\in X \st x\sim_u p \text{ for some } p\in \mathbf{P}}\\
G^s(\mathbf{P}) &= \set{(x,y)\in X^u(\mathbf{P})\times X^u(\mathbf{P}) \st x\sim_s y}.
\end{align*}
The groupoid $G^s(\mathbf{P})$ is called the {\bf stable groupoid}. Similarly, we can define
\begin{align*}
X^s(\mathbf{P}) &= \set{x\in X \st x\sim_s p \text{ for some } p\in \mathbf{P}}\\
G^u(\mathbf{P}) &= \set{(x,y)\in X^s(\mathbf{P})\times X^s(\mathbf{P}) \st x\sim_u y},
\end{align*}
and $G^u(\mathbf{P})$ is called the {\bf unstable groupoid}. There are topologies on $G^s(\mathbf{P})$ and $G^u(\mathbf{P})$ making them into 2nd countable, locally compact Hausdorff, \'etale groupoids \cite{PutSpi}. As equivalence relations, both are principal and their unit spaces are homeomorphic to $X^u(\mathbf{P})$ and $X^s(\mathbf{P})$, respectively. The corresponding groupoid $C^*$-algebras, $C^*_r(G^s(\mathbf{P}))$ and $C^*_r(G^u(\mathbf{P}))$, are called the {\bf stable algebra} and {\bf unstable algebra}, respectively.

\subsection{Wieler Solenoids}

Irreducible Smale spaces with totally disconnected stable sets are classified by Wieler in \cite{Wie} as certain inverse limit systems called solenoids. Thus, a Smale space with totally disconnected stable sets is called a {\bf Wieler solenoid}. Let $Y$ be a flat manifold and $g:Y\to Y$ an expansive self-cover. We construct a Wieler solenoid from $Y$ and $g$ by defining
\begin{align*}
X &= \varprojlim(Y,g) =\set{(y_0,y_1,y_2,\dots) \st y_i\in Y \text{ and } g(y_{i+1})=y_i},\\
\phi(y_0,y_1,y_2,\dots) &= (g(y_0),g(y_1),g(y_2),\dots) = (g(y_0),y_0,y_1,\dots).
\end{align*}
By Theorem 1 of \cite{Shub}, $Y$ always has a $g$-fixed point, which we denote by $p$. For example, for $g$ as in Theorem \ref{thm:ES}, we can take $p$ to be the equivalence class of the origin in the quotient $Y=\R^d/\pi$. The point $P=(p,p,p,\dots)\in X$ is a $\phi$-fixed point. Since $(X,\phi)$ has totally disconnected stable sets, the unstable groupoid $G^u(P)$ is ample. The unstable sets, however, are homeomorphic to $\R^d$ by Proposition 2.10 of \cite{reu} and so $G^s(P)$ is not ample.

The following two results relate the unstable groupoid of such a Wieler solenoid to the odometers considered in the rest of the paper.

\begin{prop}[\cite{reu} Lemma 3.6]
\label{prop:moritagroupoid}
Let $G^u(P)$ be the unstable groupoid for the Wieler solenoid corresponding to a flat manifold $Y$ and an expansive self-cover $g:Y\to Y$. Then $G^u(P)$ is Morita equivalent to the odometer $\Omega\rtimes \pi_1(Y)$.
\end{prop}

\begin{cor}
For each $d\geq 4$, there is a Smale space $(X,\phi)$ with $\dim (X) = d$ such that the unstable groupoid $G^u(P)$ is a counterexample to the HK-conjecture.
\end{cor}

Since $G^s(P)$ is not ample, it does not fit into the scope of the HK-conjecture. Nevertheless, there is interest in understanding when a similar statement holds for non-ample groupoids (see for example \cite{PY4}). The following theorem from \cite{reu} allows us to compute the relevant invariants.

\begin{prop}[\cite{reu} Theorem 3.2]
Let $G^s(P)$ be the stable groupoid for the Wieler solenoid corresponding to a flat manifold $Y$ and an expanding self-cover $g:Y\to Y$. Let $\operatorname{tr}^H$ (resp. $\operatorname{tr}^K$) be the transfer map in cohomology (resp. $K$-theory). Then
\begin{align*}
H_*(G^s(P)) &\iso \varinjlim(H^*(Y),\operatorname{tr}^H),\\
K_*(C^*_r(G^s(P))) &\iso \varinjlim(K^*(Y),\operatorname{tr}^K).
\end{align*}
\end{prop}

%Combining this with our main results, we obtain the following.

\begin{cor}
For each $d\geq 4$ there exists a Smale space $(X,\phi)$ with $\dim(X)=d$ such that
\[
K_*(C^*_r(G^s(P))) \not\iso \bigoplus_i H_{2i+*}(G^s(P)).
\]
\begin{proof}
Take $Y$ to be a flat manifold of dimension $d$ as in Theorem \ref{thm:largedim}, so that $T(K^*(Y)) \not\iso T(\bigoplus_i H^{2i+*}(Y))$. Take the expansive $n$-fold self-cover $g: Y\to Y$ as in Theorem \ref{thm:ES} so that for each torsion element $\gamma\in T(H^*(Y))$ we have $n\gamma=\gamma$. Then the associated odometer $\G$ is a counterexample to the HK-conjecture by Theorem \ref{thm:counterexample}. Since the cohomology transfer map satisfies $\operatorname{tr}^H\circ g^* = \times n$, we know that $\operatorname{tr}^H|_{T(H^*(Y))}:T(H^*(Y))\to T(H^*(Y))$ is an isomorphism. %Explicitly, for $\gamma \in T(H^*(Y))$, $\operatorname{tr}^H\circ g^*(\gamma) = n\gamma =\gamma$ so that $\operatorname{tr}^H$ is surjective. Since $T(H^*(Y))$ is a finite group, $\operatorname{tr}^H$ is also injective.
Thus, by Lemma \ref{lem:torlim}, we have $T\parens{H_{*}(G^s(P))} \iso T\parens{\varinjlim(H^*(Y), \operatorname{tr}^H)} \iso T(H^*(Y)),$ where the first isomorphism is by the previous proposition. Also by Lemma \ref{lem:torlim}, $T(K_*(C^*_r(G^s(P))))$ is a subgroup of $T(K^*(Y))$. By Proposition \ref{prop:Ktorsion}, $|T(K^*(Y))|\leq |T(\bigoplus_i H^{2i+*}(Y))|$, so that $T(K_*(C^*_r(G^s(P)))) \not\iso T(\bigoplus_i H^{2i+*}(Y))$ and hence the desired result.
\end{proof}
\end{cor}

However, an HK-type statement does hold for some low-dimensional examples.% such as $Y=S^1$.

\begin{example}
Consider $Y=S^1$ and $g:S^1\to S^1$ the $m$-fold cover induced by multiplication by $m$ on $\R$. The Wieler solenoid $X=\varprojlim(Y,g)$ is known as the $m$-solenoid. In \cite{Put}, Putnam shows that the stable and unstable algebras for the $m$-solenoid are both isomorphic to the stabilized Bunce--Deddens algebra of type $m^{\infty}$. Thus, $K_0(C^*_r(G^s(P)))\iso \Z\bracks{\frac{1}{m}}$ and $K_1(C^*_r(G^s(P)))\iso \Z$ (see \cite{Blackadar}). Theorems 3.10 and 3.11 of \cite{reu} give $H_0(G^s(P)) \iso \Z\bracks{\frac{1}{m}}$ and $H_1(G^s(P)) \iso \Z,$ with all other homology groups trivial. This gives the HK-type statement for the $m$-solenoid.

We remark that this also follows from Proposition 3.2 of \cite{MR4052913}, using the fact that the stable and unstable algebras are isomorphic and that the two theorems from \cite{reu} also give $H_n(G^s(P))\iso H_n(G^u(P))$ for $n=0,1$.
\end{example}

%%%%%%%%%%%%%%%%% Bibliography %%%%%%%%%%%%%%%%%
\bibliographystyle{amsplain}
\bibliography{Biblio-Database}

\end{document}